\documentclass{amsart}
\usepackage{amssymb}
\usepackage{amsmath}
\usepackage{amsfonts}

\newtheorem{theorem}{Theorem}
\theoremstyle{plain}

\newtheorem{definition}[theorem]{Definition}

\newtheorem{lemma}[theorem]{Lemma}
\newtheorem{notation}[theorem]{Notation}

\newtheorem{proposition}[theorem]{Proposition}
\newtheorem{remark}[theorem]{Remark}

\numberwithin{equation}{section}
\numberwithin{theorem}{section}
\input{tcilatex}

\begin{document}
\title[Quadratic reciprocity and Gauss sum via the Weil representation]{%
Quadratic reciprocity and the sign of the Gauss sum via the finite Weil
representation}
\author{Shamgar Gurevich}
\address{Department of Mathematics, University of California, Berkeley, CA
94720, USA.}
\email{shamgar@math.berkeley.edu}
\author{Ronny Hadani}
\address{Department of Mathematics, University of Chicago, IL 60637, USA.}
\email{hadani@math.uchicago.edu}
\author{Roger Howe}
\address{Department of Mathematics, Yale University, New Haven, CT
06520-8283, USA.}
\email{howe@math.yale.edu}
\date{January 1, 2008.}
\thanks{\copyright \ Copyright by S. Gurevich, R. Hadani and R. Howe, January
1, 2008. All rights reserved.}

\begin{abstract}
We give new proofs of two basic results in number theory: The law of
quadratic reciprocity and the sign of the Gauss sum. We show that these
results are encoded in the relation between the discrete Fourier transform
and the action of the Weyl element in the Weil representation modulo $p,q$
and $pq.$
\end{abstract}

\maketitle

\section{Introduction}

Two basic results due to Gauss are the quadratic reciprocity law and the
sign of the Gauss sum \cite{IR}. The first concerns the identity%
\begin{equation}
\QOVERD( ) {p}{q}\QOVERD( ) {q}{p}=(-1)^{\frac{p-1}{2}\frac{q-1}{2}},
\label{QRF}
\end{equation}%
where $p,q$ are two distinct odd prime numbers and $\QOVERD( ) {\cdot }{p}$
(respectively $\QOVERD( ) {\cdot }{q}$) is the Legendre symbol modulo $p$
(respectively $q$), i.e., $\QOVERD( ) {x}{p}=1$ if $x$ is a square modulo $p$
and $-1$ otherwise. The latter result asserts that%
\begin{equation}
G_{p}=\dsum \limits_{x\in \mathbb{F}_{p}}e^{\frac{2\pi i}{p}x^{2}}=\left \{ 
\begin{array}{cc}
\sqrt{p}, & p\equiv 1\text{ }(\func{mod}\text{ }4), \\ 
i\sqrt{p}, & p\equiv 3\text{ }(\func{mod}\text{ }4).%
\end{array}%
\right.  \label{GSF}
\end{equation}

In fact, it is easy to show that $G_{p}^{2}=\QOVERD( ) {-1}{p}\cdot p.$
Hence, the problem is to determine the exact sign in the evaluation of $%
G_{p}.$

In this work we will explain how these results follow from the
proportionality relations 
\begin{equation*}
F_{n}=C_{n}\cdot \rho _{n}%
\begin{pmatrix}
0 & -1 \\ 
1 & 0%
\end{pmatrix}%
,\text{ for }n=p,q\text{ and }pq,
\end{equation*}%
where $F_{n}$ is the discrete Fourier transform, $\rho _{n}$ is the Weil
representation of the group $SL_{2}(\mathbb{%
\mathbb{Z}
}/n\mathbb{%
\mathbb{Z}
}),$ both acting on the Hilbert space $%
\mathbb{C}
(\mathbb{%
\mathbb{Z}
}/n\mathbb{%
\mathbb{Z}
})$ of complex valued functions on the finite ring $\mathbb{%
\mathbb{Z}
}/n\mathbb{%
\mathbb{Z}
},$ and $C_{n}$ is the proportionality constant. More specifically, the law
of quadratic reciprocity follows from basic properties of the Weil
representation and group theoretic considerations, while the calculation of
the sign of the Gauss sum is a bit more delicate and it uses a formula for
the character of the Weil representation. The fact that the discrete Fourier
transform can be normalized so that it becomes a part of a representation
plays a crucial role in our proof of both statements.

In his seminal work \cite{W}, Andr\'{e} Weil recast several known proofs of
the law of quadratic reciprocity in terms of the Weil representation of some
cover of the group $SL_{2}(\mathbb{A}_{%
\mathbb{Q}
})$ where $\mathbb{A}_{%
\mathbb{Q}
}$ denotes the adele ring of $%
\mathbb{Q}
.$ The main contribution of this short note is showing that quadratic
reciprocity already follows from the Weil representation over finite rings
and, moreover, establishing a conceptual mechanism, different from that of
Weil, which produces the law of quadratic reciprocity.

\subsection{Structure of the paper}

In Section \ref{WR} we recall the Weil representation over the finite ring $%
\mathbb{Z}
/n%
\mathbb{Z}
.$ We then describe the relation between the Weil representations associated
with the rings $%
\mathbb{Z}
/n_{1}%
\mathbb{Z}
$, $%
\mathbb{Z}
/n_{2}%
\mathbb{Z}
$ and $%
\mathbb{Z}
/n_{1}n_{2}%
\mathbb{Z}
$, for $n_{1},n_{2}$ coprime. Finally, we write the formula of its character
in the case $n$ is an odd prime number. In Section \ref{DFTandW} we define
the discrete Fourier transform, compute its determinant and explain its
relation to the Weil representation. In Section \ref{quad_sec} we prove the
quadratic reciprocity law and in Section \ref{Gauss_sec} we compute the
Gauss sum. Finally, in the Appendix we supply the proofs of the main
technical claims that appear in the body of the paper.

\subsection{Acknowledgements}

The first two authors would like to thank their teacher J. Bernstein for his
interest and guidance. They would like also to acknowledge M. Nori for the
encouragement to write this paper. They thank C.P. Mok for explaining parts
from the known proofs of quadratic reciprocity and T. Schedler for helping
with various of the computations. We appreciate the time T.Y. Lam spent with
us teaching math and history. Finally, we thank M. Baruch, P. Diaconis, M.
Haiman, B. Poonen and K. Ribet for the opportunities to present this work at
the Technion, Israel, and the MSRI, RTG and number theory seminars at\
Berkeley during February 2008.

\section{The Weil representation\label{WR}}

\subsection{The Heisenberg group}

Let $(V,\omega )$ be a symplectic free module of rank $2$ over the finite
ring $\mathbb{%
\mathbb{Z}
}/n%
\mathbb{Z}
$, where $n$ is an odd number. The reader should think of $V$ as $\mathbb{%
\mathbb{Z}
}/n%
\mathbb{Z}
\times \mathbb{%
\mathbb{Z}
}/n%
\mathbb{Z}
$ equipped with the standard skew-symmetric form $\omega \left( \left(
t,w\right) ,\left( t^{\prime },w^{\prime }\right) \right) =tw^{\prime
}-wt^{\prime }$. Considering $V$ as an abelian group, it admits a
non-trivial central extension called the \textit{Heisenberg }group. \
Concretely, the group $H$ can be presented as the set $H=V\times \mathbb{%
\mathbb{Z}
}/n%
\mathbb{Z}
$ with the multiplication given by%
\begin{equation*}
(v,z)\cdot (v^{\prime },z^{\prime })=(v+v^{\prime },z+z^{\prime }+\tfrac{1}{2%
}\omega (v,v^{\prime })).
\end{equation*}

The center of $H$ is $Z=Z(H)=\left \{ (0,z):\text{ }z\in \mathbb{%
\mathbb{Z}
}/n%
\mathbb{Z}
\right \} .$ The symplectic group $Sp(V,\omega )$, which in this case is
isomorphic to $SL_{2}\left( \mathbb{%
\mathbb{Z}
}/n%
\mathbb{Z}
\right) $, acts by automorphisms of $H$ through its action on the $V$%
-coordinate.

\subsection{The Heisenberg representation\label{HR}}

One of the most important attributes of the group $H$ is that it admits a
special family of irreducible representations. The precise statement goes as
follows. Let $\psi :Z\rightarrow 
\mathbb{C}
^{\times }$ be a faithful character of the center (i.e., an imbedding of $Z$
into $%
\mathbb{C}
^{\times })$. It is not hard to show

\begin{theorem}[Stone-von Neuman]
\label{S-vN}There exists a unique (up to isomorphism) irreducible
representation $(\pi ,H,\mathcal{H)}$ with the center acting by $\psi ,$
i.e., $\pi _{|Z}=\psi \cdot Id_{\mathcal{H}}$.
\end{theorem}

The representation $\pi $ which appears in the above theorem will be called
the \textit{Heisenberg representation associated with the central character }%
$\psi $.

We denote by $\psi _{1}(z)=e^{\frac{2\pi i}{n}z}$ the standard additive
character, and for every invertible element $a\in \left( \mathbb{%
\mathbb{Z}
}/n%
\mathbb{Z}
\right) ^{\times }$ we denote $\psi _{a}\left( z\right) =e^{\frac{2\pi i}{n}%
az}$.

\subsubsection{Standard realization of the Heisenberg representation\label%
{standard_subsub}.}

The Heisenberg representation $(\pi ,H,\mathcal{H)}$ can be realized as
follows: The Hilbert space is the space $\mathcal{F}_{n}=%
\mathbb{C}
(\mathbb{%
\mathbb{Z}
}/n%
\mathbb{Z}
)$ of complex valued functions on the finite field, with the standard
Hermitian product. The action $\pi $ is given by

\begin{itemize}
\item $\pi (t,0)[f]\left( x\right) =f\left( x+t\right) ;$

\item $\pi (0,w)[f]\left( x\right) =\psi \left( wx\right) f\left( x\right) ;$

\item $\pi (z)[f]\left( x\right) =\psi \left( z\right) f\left( x\right) ,$ $%
z\in Z.$
\end{itemize}

Here we are using $t$ to indicate the first coordinate of a typical element $%
v\in V\simeq \mathbb{%
\mathbb{Z}
}/n%
\mathbb{Z}
\times \mathbb{%
\mathbb{Z}
}/n%
\mathbb{Z}
$, and $w$ to indicate the second coordinate.

We call this explicit realization the \textit{standard realization}.

\subsection{The Weil representation\label{Wrep_sub}}

A direct consequence of Theorem \ref{S-vN} is the existence of a projective
representation $\widetilde{\rho }_{n}:SL_{2}(\mathbb{%
\mathbb{Z}
}/n%
\mathbb{Z}
)\rightarrow PGL(\mathcal{H)}$. The construction of $\widetilde{\rho }_{n}$
out of the Heisenberg representation $\pi $ is due to Weil \cite{W} and it
goes as follows: Considering the Heisenberg representation $\pi $ and an
element $g\in SL_{2}(\mathbb{%
\mathbb{Z}
}/n%
\mathbb{Z}
)$, one can define a new representation $\pi ^{g}$ acting on the same
Hilbert space via $\pi ^{g}\left( h\right) =\pi \left( g\left( h\right)
\right) $. Clearly both $\pi $ and $\pi ^{g}$ have the same central
character $\psi ,$ hence, by Theorem \ref{S-vN}, they are isomorphic. Since
the space $\mathsf{Hom}_{H}(\pi ,\pi ^{g})$ is one-dimensional, choosing for
every $g\in SL_{2}(\mathbb{%
\mathbb{Z}
}/n%
\mathbb{Z}
)$ a non-zero representative $\widetilde{\rho }_{n}(g)\in \mathsf{Hom}%
_{H}(\pi ,\pi ^{g})$ gives the required projective Weil representation. In
more concrete terms, the projective representation $\widetilde{\rho }$ is
characterized by the formula 
\begin{equation}
\widetilde{\rho }_{n}\left( g\right) \pi \left( h\right) \widetilde{\rho }%
_{n}\left( g\right) ^{-1}=\pi \left( g\left( h\right) \right) ,
\label{Egorov}
\end{equation}%
for every $g\in SL_{2}(\mathbb{%
\mathbb{Z}
}/n%
\mathbb{Z}
)$ and $h\in H$. A more delicate statement is that there exists a lifting of 
$\widetilde{\rho }_{n}$ into a linear representation, this is the content of
the following theorem

\begin{theorem}
\label{linearization_thm} The projective Weil representation $\widetilde{%
\rho }_{n}$ can be linearized into an honest representation 
\begin{equation*}
\rho _{n}:SL_{2}(\mathbb{%
\mathbb{Z}
}/n%
\mathbb{Z}
)\longrightarrow GL(\mathcal{H)},
\end{equation*}%
that satisfies equation (\ref{Egorov}).
\end{theorem}

The existence of a linearization $\rho _{n}$ follows from a known fact \cite%
{B} that any projective representation of $SL_{2}(\mathbb{%
\mathbb{Z}
}/n%
\mathbb{Z}
)$ can be linearized (in case $n=p$ is a prime number, see also \cite{GH1,
GH2} for an explicit construction of a canonical linearization).

Clearly, any two linearizations differ by a character $\chi $ of $SL_{2}(%
\mathbb{%
\mathbb{Z}
}/n%
\mathbb{Z}
)$. In addition, we have

\begin{proposition}
\label{order_prop}Let $\chi $ be a character of the group $SL_{2}(\mathbb{%
\mathbb{Z}
}/n%
\mathbb{Z}
)$ then $\chi ^{n}=1$.
\end{proposition}

For a proof, see Appendix \ref{proofs_sec}.

\begin{remark}
\label{perfect_rmk}In the case when $n$ is not divisible by $3$ the group $%
SL_{2}(\mathbb{%
\mathbb{Z}
}/n%
\mathbb{Z}
)$ is perfect therefore the representation $\rho _{n}$ is unique. The
perfectness of $SL_{2}(\mathbb{%
\mathbb{Z}
}/n%
\mathbb{Z}
)$ can be proved as follows: Let $pr:SL_{2}(\mathbb{%
\mathbb{Z}
})\rightarrow SL_{2}(\mathbb{%
\mathbb{Z}
}/n%
\mathbb{Z}
)$ denote the canonical projection. Given a character $\chi :SL_{2}(\mathbb{%
\mathbb{Z}
}/n%
\mathbb{Z}
)\rightarrow 
\mathbb{C}
^{\times }$, its pull-back $\chi \circ pr$ satisfies $\left( \chi \circ
pr\right) ^{12}=1$ since the group of characters of $SL_{2}(\mathbb{%
\mathbb{Z}
})$ is isomorphic to $%
\mathbb{Z}
/12%
\mathbb{Z}
$ \cite{CM}. Since $pr$ is surjective this implies that $\chi ^{12}=1$. This
combined with the facts that $\chi ^{n}=1$ (Proposition \ref{order_prop})
and $\gcd (12,n)=1$ implies that $\chi =1.$
\end{remark}

\begin{notation}
The Weil representation $\rho _{n}$ depends on the central character $\psi $%
, hence, sometimes we will write $\rho _{n}\left[ \psi \right] $ to
emphasize this point. We will denote by $\rho _{n}$ the Weil representation
associated with the standard character $\psi _{1}$.
\end{notation}

Let $n_{1},n_{2}$ be coprime odd integers. Consider the natural homomorphism 
\begin{equation*}
\mathbb{Z}
/n_{1}n_{2}%
\mathbb{Z}
\overset{\simeq }{\rightarrow }%
\mathbb{Z}
/n_{1}%
\mathbb{Z}
\times 
\mathbb{Z}
/n_{2}%
\mathbb{Z}
,
\end{equation*}
which, by the Chinese reminder theorem, is an isomorphism. This isomorphism
induces an isomorphism of Hilbert spaces $\mathcal{F}_{n_{1}n_{2}}\overset{%
\simeq }{\rightarrow }\mathcal{F}_{n_{1}}\otimes \mathcal{F}_{n_{2}}$ and an
isomorphism of groups $SL_{2}(\mathbb{%
\mathbb{Z}
}/n_{1}n_{2}%
\mathbb{Z}
)\overset{\simeq }{\rightarrow }SL_{2}(\mathbb{%
\mathbb{Z}
}/n_{1}%
\mathbb{Z}
)\times SL_{2}(\mathbb{%
\mathbb{Z}
}/n_{2}%
\mathbb{Z}
)$. In addition, the character $\psi _{1}\otimes \psi _{1}$ of $%
\mathbb{Z}
/n_{1}%
\mathbb{Z}
\times 
\mathbb{Z}
/n_{2}%
\mathbb{Z}
$ transforms to the character $\psi _{n_{1}+n_{2}}$ of $%
\mathbb{Z}
/n_{1}n_{2}%
\mathbb{Z}
$. Under these identifications it is not difficult to show

\begin{proposition}
\label{tensor_prop}The representations $\rho _{n_{1}}\otimes \rho _{n_{2}}$
and $\rho _{n_{1}n_{2}}\left[ \psi _{n_{1}+n_{2}}\right] $, realized on the
Hilbert spaces $\mathcal{F}_{n_{1}}\otimes \mathcal{F}_{n_{2}}$ and $%
\mathcal{F}_{n_{1}+n_{2}}$, coincide as projective representations of $%
SL_{2}(\mathbb{%
\mathbb{Z}
}/n_{1}n_{2}%
\mathbb{Z}
)$.
\end{proposition}

\begin{remark}
The element $n_{1}+n_{2}\in \left( \mathbb{%
\mathbb{Z}
}/n_{1}n_{2}%
\mathbb{Z}
\right) ^{\times }$, hence, the character $\psi _{n_{1}+n_{2}}$ is faithful.
\end{remark}

\begin{remark}
In the case $n_{1},n_{2}$ are in addition not divisible by $3$, Proposition %
\ref{tensor_prop} and Remark \ref{perfect_rmk} imply that the
representations $\rho _{n_{1}n_{2}}\left[ \psi _{n_{1}+n_{2}}\right] $ and $%
\rho _{n_{1}}\left[ \psi _{1}\right] \otimes \rho _{n_{2}}\left[ \psi _{1}%
\right] $ coincide.
\end{remark}

\subsection{The character of the Weil representation}

In the case $n=p$ is a prime number, the absolute value of the character $%
ch_{\rho }:$ $SL_{2}(\mathbb{F}_{p})\rightarrow 
\mathbb{C}
$ of the Weil representation $\rho =\rho _{p}$ was described in \cite{H},
but the phases have been made explicit only recently in \cite{GH1}. The
following formula is taken from \cite{GH1}:

\begin{equation}
ch_{\rho }(g)=\QOVERD( ) {-\det \left( \kappa \left( g\right) +I\right) }{p},
\label{Ch}
\end{equation}%
for every $g\in SL_{2}(\mathbb{F}_{p})$ such that $g-I$ is invertible, where 
$\QOVERD( ) {\cdot }{p}$ is the Legendre symbol modulo $p$ and $\kappa
\left( g\right) =\frac{g+I}{g-I}$ is the Cayley transform. Using the
identity $\kappa \left( g\right) +I=\frac{2g}{g-I}$ one can write (\ref{Ch})
in the simpler form%
\begin{equation*}
ch_{\rho }(g)=\QOVERD( ) {-\det \left( g-I\right) }{p}.
\end{equation*}

\begin{remark}
Sketch of the proof of (\ref{Ch}) (see details in \cite{GH1}). First
observation is that the Heisenberg representation $\pi $ and the Weil
representation $\rho =\rho _{p}$ combine to give a representation of the
semi-direct product $\tau =\rho \ltimes \pi :SL_{2}\left( \mathbb{F}%
_{p}\right) \ltimes H\rightarrow GL\left( \mathcal{H}\right) $. Second
observation\ is that the character of $\tau $, $ch_{\tau }:SL_{2}\left( 
\mathbb{F}_{p}\right) \ltimes H\rightarrow 
\mathbb{C}
$ satisfies the following multiplicativity property%
\begin{equation}
ch_{\tau }\left( g_{1}\cdot g_{2}\right) =\frac{1}{\dim \mathcal{H}}ch_{\tau
}\left( g_{1}\right) \ast ch_{\tau }\left( g_{2}\right) ,  \label{conv_eq}
\end{equation}%
where $ch_{\tau }\left( g\right) $ denotes the function on $H$ given by $%
ch_{\tau }\left( g,-\right) $ and the $\ast $ operation denotes convolution
with respect to the Heisenberg group action. Now, one can easily show that 
\begin{equation*}
ch_{\tau }\left( g\right) \left( v,z\right) =\mu _{g}\cdot \psi \left( 
\tfrac{1}{4}\omega \left( \kappa \left( g\right) v,v\right) +z\right) ,
\end{equation*}%
for $g\in $ $SL_{2}\left( \mathbb{F}_{p}\right) $ with $g-I$ invertible, and
for some $\mu _{g}\in 
\mathbb{C}
$. Finally, a direct calculation reveals that $\mu _{g}$ must equal $%
\QOVERD( ) {-\det \left( \kappa \left( g\right) +I\right) }{p}$ for (\ref%
{conv_eq}) to hold. Restricting $ch_{\tau }$ to $SL_{2}\left( \mathbb{F}%
_{p}\right) \subset SL_{2}\left( \mathbb{F}_{p}\right) \ltimes H$ we obtain (%
\ref{Ch}).
\end{remark}

\section{The discrete Fourier transform\label{DFTandW}}

\subsection{ The discrete Fourier transform (DFT)}

Given an additive character $\psi $, there is a DFT operator $F_{n}\left[
\psi \right] $ acting on the Hilbert space $\mathcal{F}_{n}\mathcal{=%
\mathbb{C}
}\left( \mathbb{Z}/n\mathbb{Z}\right) $ by the formula\footnote{%
Usually, the DFT operator appears in its normalized form $\Phi _{n}=\frac{1}{%
\sqrt{n}}F_{n}$, which makes it a unitary operator. However, the
non-normalized form is better suited to our purposes.} 
\begin{equation*}
F_{n}\left[ \psi \right] \left( f\right) \left( y\right) =\sum \limits_{x\in 
\mathbb{%
\mathbb{Z}
}/n\mathbb{%
\mathbb{Z}
}}\psi \left( yx\right) f\left( x\right) .
\end{equation*}%
It is easy to show that $F_{n}\left[ \psi \right] ^{2}(f)\left( y\right)
=n\cdot f(-y)$ and hence $F_{n}\left[ \psi \right] ^{4}=n^{2}\cdot Id.$

\begin{notation}
We will denote by $F_{n}$ the DFT operator associated with the standard
character $\psi _{1}$.
\end{notation}

It is not difficult to show that under the isomorphism $\mathcal{F}%
_{n_{1}n_{2}}\overset{\simeq }{\rightarrow }\mathcal{F}_{n_{1}}\otimes 
\mathcal{F}_{n_{2}}$ we have

\begin{proposition}
\label{tensor2_prop}The operators $F_{n_{1}}\otimes F_{n_{2}}$ and $%
F_{n_{1}n_{2}}\left[ \psi _{n_{1}+n_{2}}\right] $ coincide.
\end{proposition}

The main technical statement that we will require concerns the explicit
evaluation of the determinant of the DFT operator which is associated with
the standard character $\psi _{1}$.

\begin{proposition}
\label{detF}For any odd natural number $n$ we have 
\begin{equation}
\det \left( F_{n}\right) =i^{\frac{n(n-1)}{2}}n^{\frac{n}{2}}.
\label{det(F)}
\end{equation}
\end{proposition}

For a proof, see Appendix \ref{proofs_sec}.

\subsection{Relation between the DFT and the Weil representation}

We will show that for an odd $n$ the operator $F_{n}$ is proportional to the
operator $\rho _{n}\left( \mathrm{w}\right) $ in the Weil representation,
where $\mathrm{w}\in SL_{2}\left( \mathbb{%
\mathbb{Z}
}/n%
\mathbb{Z}
\right) $ is the Weyl element 
\begin{equation*}
\mathrm{w}=%
\begin{pmatrix}
0 & 1 \\ 
-1 & 0%
\end{pmatrix}%
.
\end{equation*}

\begin{theorem}
\label{DFT_lemma}The operator $\rho _{n}\left( \mathrm{w}\right) $ does not
depend on the choice of linearization $\rho _{n}$. Moreover, 
\begin{equation}
F_{n}=C_{n}\cdot \rho _{n}\left( \mathrm{w}\right) ,  \label{poportion_eq}
\end{equation}%
where $C_{n}=i^{\frac{n-1}{2}}\sqrt{n}$.
\end{theorem}

For a proof see Appendix \ref{proofs_sec}.

\section{The quadratic reciprocity law\label{quad_sec}}

We are ready to prove the quadratic reciprocity law. For a natural number $n$%
, we define the Gauss sums \cite{IR} 
\begin{equation*}
G_{n}\left( a\right) =\dsum \limits_{x\in \mathbb{%
\mathbb{Z}
}/n%
\mathbb{Z}
}e^{\frac{2\pi i}{n}ax^{2}},
\end{equation*}%
where $a\in \left( \mathbb{%
\mathbb{Z}
}/n%
\mathbb{Z}
\right) ^{\times }$ and we denote $G_{n}=G_{n}\left( 1\right) $.

Consider two distinct odd prime numbers $p,q$.

\begin{lemma}
\label{ident_lemma}We have 
\begin{equation*}
\QOVERD( ) {p}{q}\QOVERD( ) {q}{p}=\frac{G_{p}\cdot G_{q}}{G_{pq}}.
\end{equation*}
\end{lemma}

For a proof see Appendix \ref{proofs_sec}.

The Gauss sum is related to the DFT by 
\begin{equation}
G_{n}=Tr(F_{n}).  \label{G-F}
\end{equation}

Hence we can write%
\begin{eqnarray}
\QOVERD( ) {p}{q}\QOVERD( ) {q}{p} &=&\frac{Tr(F_{p})\cdot Tr(F_{q})}{%
Tr(F_{pq})}=\frac{C_{p}\cdot C_{q}}{C_{pq}}\cdot \frac{Tr(\rho _{p}(\mathrm{w%
}))\cdot Tr(\rho _{q}(\mathrm{w}))}{Tr(\rho _{pq}(\mathrm{w}))}  \notag \\
&=&\frac{C_{p}\cdot C_{q}}{C_{pq}}=(-1)^{\frac{p-1}{2}\cdot \frac{q-1}{2}},
\label{QR_eq}
\end{eqnarray}%
where in the first equality we used Equation (\ref{G-F}) and Lemma \ref%
{ident_lemma}, in the second and forth equalities we used Theorem \ref%
{DFT_lemma}. Finally, the third equality follows from

\begin{proposition}
\label{trace_prop}We have 
\begin{equation*}
Tr(\rho _{pq}(\mathrm{w}))=Tr(\rho _{p}(\mathrm{w}))\cdot Tr(\rho _{q}(%
\mathrm{w})).
\end{equation*}
\end{proposition}

For a proof, see Appendix \ref{proofs_sec}.

This completes our proof of Equation (\ref{QRF}) - the Quadratic Reciprocity
Law.

\subsection{Quadratic reciprocity law for the Jacobi symbol}

For an odd number $n\in 
\mathbb{N}
$, let $\QOVERD( ) {\cdot }{n}$ denote the Jacobi symbol of the
multiplicative group $\left( 
\mathbb{Z}
/n%
\mathbb{Z}
\right) ^{\times }$, which can be characterized \cite{BEW} by the condition 
\begin{equation}
G_{n}\left( a\right) =\left( \frac{a}{n}\right) G_{n}\left( 1\right) ,
\label{equiv_eq}
\end{equation}%
for every $a\in \left( 
\mathbb{Z}
/n%
\mathbb{Z}
\right) ^{\times }$.

\begin{remark}
The Jacobi symbol admits the following explicit description: When $n=p$ is
an odd prime number, the Jacobi symbol coincides with the Legendre character
(Lemma \ref{tech_lemma}). If $n=p_{1}^{k_{1}}\cdot ..\cdot p_{l}^{k_{l}}$ is
the decomposition of $n$ into a product of prime numbers then it can be
shown that for $a\in \left( 
\mathbb{Z}
/n%
\mathbb{Z}
\right) ^{\times }$, 
\begin{equation*}
\left( \frac{a}{n}\right) =\left( \frac{a}{p_{1}}\right) ^{k_{1}}\cdot
..\cdot \left( \frac{a}{p_{l}}\right) ^{k_{l}}\text{.}
\end{equation*}%
In particular, this implies that the Jacobi symbol is a character of the
multiplicative group $\left( 
\mathbb{Z}
/n%
\mathbb{Z}
\right) ^{\times }$ and it takes the values $\pm 1$.
\end{remark}

The quadratic reciprocity law can be formulated in terms of the Jacobi
symbol, for any two coprime odd numbers $n_{1},n_{2}$; the general law is 
\begin{equation}
\QOVERD( ) {n_{1}}{n_{2}}\QOVERD( ) {n_{2}}{n_{1}}=(-1)^{\frac{n_{1}-1}{2}%
\cdot \frac{n_{2}-1}{2}}.  \label{QR_gen}
\end{equation}

The proof we just described for the quadratic reciprocity law gives also the
more general identity (\ref{QR_gen}) without changes: Using Equation (\ref%
{equiv_eq}), it is not hard to realize that the statement of Lemma \ref%
{ident_lemma} can be formulated more generally as 
\begin{equation*}
\QOVERD( ) {n_{1}}{n_{2}}\QOVERD( ) {n_{2}}{n_{1}}=\frac{G_{n_{1}}\cdot
G_{n_{2}}}{G_{n_{1}n_{2}}}.
\end{equation*}

Then applying the same derivation as in (\ref{QR_eq}) one obtains 
\begin{equation*}
\QOVERD( ) {n_{1}}{n_{2}}\QOVERD( ) {n_{2}}{n_{1}}=\frac{C_{n_{1}}\cdot
C_{n_{2}}}{C_{n_{1}n_{2}}}=(-1)^{\frac{n_{1}-1}{2}\cdot \frac{n_{2}-1}{2}}.
\end{equation*}

\subsection{Alternative interpretation}

A slightly more transparent interpretation of the above argument proceeds as
follows: 
\begin{equation*}
1=\frac{Tr\left( F_{n_{1}}\right) Tr\left( F_{n_{2}}\right) }{Tr\left(
F_{n_{1}n_{2}}\left[ \psi _{n_{1}+n_{2}}\right] \right) }=\frac{%
C_{n_{1}}C_{n_{2}}}{C_{n_{1}n_{2}}\left[ \psi _{n_{1}+n_{2}}\right] }\cdot 
\frac{Tr\left( \rho _{n_{1}}\left( \mathrm{w}\right) \right) Tr\left( \rho
_{n_{2}}\left( \mathrm{w}\right) \right) }{Tr\left( \rho _{n_{1}n_{2}}\left[
\psi _{n_{1}+n_{2}}\right] \left( \mathrm{w}\right) \right) }\text{,}
\end{equation*}%
where the first equality follows from Proposition \ref{tensor2_prop} and the
second equality appears by substituting 
\begin{eqnarray*}
F_{n_{1}} &=&C_{n_{1}}\rho _{n_{1}}\left( \mathrm{w}\right) , \\
F_{n_{2}} &=&C_{n_{2}}\rho _{n_{2}}\left( \mathrm{w}\right) , \\
F_{n_{1}n_{2}}\left[ \psi _{n_{1}+n_{2}}\right] &=&C_{n_{1}n_{2}}\left[ \psi
_{n_{1}+n_{2}}\right] \rho _{n_{1}n_{2}}\left[ \psi _{n_{1}+n_{2}}\right]
\left( \mathrm{w}\right) .
\end{eqnarray*}

Now, by Proposition \ref{tensor_prop}, we have 
\begin{equation*}
\frac{Tr\left( \rho _{n_{1}}\left( \mathrm{w}\right) \right) Tr\left( \rho
_{n_{2}}\left( \mathrm{w}\right) \right) }{Tr\left( \rho _{n_{1}n_{2}}\left[
\psi _{n_{1}+n_{2}}\right] \left( \mathrm{w}\right) \right) }=1\text{.}
\end{equation*}

Hence, the quadratic reciprocity law follows from the equivariance property
of the proportion constant $C$:

\begin{theorem}
\label{equivariance_thm}Let $n\in 
\mathbb{N}
$ be an odd number. We have 
\begin{equation*}
C_{n}\left[ \psi _{a}\right] =\left( \frac{a}{n}\right) C_{n}\left[ \psi _{1}%
\right] \text{,}
\end{equation*}%
for every $a\in \left( 
\mathbb{Z}
/n%
\mathbb{Z}
\right) ^{\times }$.
\end{theorem}

\begin{remark}
In our approach, the statement of Theorem \ref{equivariance_thm} follows,
indirectly, from the proof of the quadratic reciprocity law.
\end{remark}

\section{The sign of Gauss sum\label{Gauss_sec}}

The exact evaluation of $G_{p}$ for an odd prime $p$ uses an additional fact
about the Weil representation, i.e., the evaluation 
\begin{equation*}
Tr(\rho _{p}(\mathrm{w}))=ch_{\rho _{p}}(\mathrm{w})=\QOVERD( ) {-2}{p},
\end{equation*}
which follows directly from formula (\ref{Ch}). The explicit evaluation of
the Legendre symbol at $-2$ is a simple and well-known computation (see \cite%
{S}) which gives $\QOVERD( ) {-2}{p}=\QOVERD( ) {-1}{p}\QOVERD( )
{2}{p}=(-1)^{\frac{p-1}{2}}(-1)^{\frac{p^{2}-1}{8}}.$ Now using (\ref{G-F})
we conclude that%
\begin{eqnarray*}
G_{p} &=&Tr(F_{p})=C_{p}\cdot Tr(\rho _{p}(\mathrm{w})) \\
&=&i^{\frac{p-1}{2}}\sqrt{p}\cdot \QOVERD( ) {-2}{p}=\left \{ 
\begin{array}{cc}
\sqrt{p}, & p\equiv 1\text{ }(\func{mod}\text{ }4), \\ 
i\sqrt{p}, & p\equiv 3\text{ }(\func{mod}\text{ }4).%
\end{array}%
\right.
\end{eqnarray*}

This completes our proof of equation (\ref{GSF}) - the Sign of Gauss Sum.

\appendix

\section{Proof of statements\label{proofs_sec}}

\subsection{Proof of Proposition \protect \ref{order_prop}}

Let $\chi $ be a character of $SL_{2}(\mathbb{%
\mathbb{Z}
}/n%
\mathbb{Z}
)$. The condition $\chi ^{n}=1$ follows from the basic fact that the group $%
SL_{2}(\mathbb{%
\mathbb{Z}
}/n%
\mathbb{Z}
)$ is generated by the unipotent elements \cite{CM}: 
\begin{equation*}
u_{+}=%
\begin{pmatrix}
1 & 1 \\ 
0 & 1%
\end{pmatrix}%
,\text{ \ }u_{-}=%
\begin{pmatrix}
1 & 0 \\ 
1 & 1%
\end{pmatrix}%
,
\end{equation*}%
which satisfy $u_{+}^{n}=u_{-}^{n}=Id$.

\subsection{Proof of Proposition \protect \ref{detF}}

Consider the matrix of $F_{n}$. If we write $\psi $ in the form $\psi \left(
x\right) =\zeta ^{x}$, $x=0,..,n-1$ and $\zeta =e^{\frac{2\pi i}{n}}$, the
matrix of $F_{n}$ takes the form $(\zeta ^{yx}:$ $\ x,y\in \left \{
0,1,..,n-1\right \} )$, hence it is a \textit{Vandermonde matrix.} Applying
the standard formula for the determinant of a Vandermonde matrix (see \cite%
{FH}), we get%
\begin{eqnarray}
\det (F_{n}) &=&\tprod \limits_{0\leq x<y\leq n-1}\left( \zeta ^{y}-\zeta
^{x}\right) =\tprod \limits_{0\leq x<y\leq n-1}(\psi (y)-\psi (x))
\label{CalcDet} \\
&=&\tprod \limits_{0\leq x<y\leq n-1}\psi (\tfrac{x+y}{2})\tprod
\limits_{0\leq x<y\leq n-1}\left( \psi (\tfrac{y-x}{2})-\psi (\tfrac{x-y}{2}%
)\right)  \notag \\
&=&\psi (0)\cdot i^{\frac{n(n-1)}{2}}\cdot 2^{\frac{n(n-1)}{2}}\tprod
\limits_{j=1}^{n-1}(\sin (\tfrac{\pi \cdot j}{n}))^{n-j},  \notag
\end{eqnarray}%
where the equality $\tprod \limits_{0\leq x<y\leq n-1}\psi (\tfrac{x+y}{2}%
)=\psi (0)$ follows from the fact that 
\begin{equation*}
\tsum \limits_{0\leq x<y\leq n-1}(x+y)=\frac{1}{2}\{ \tsum \limits_{x,y\in 
\mathbb{Z}
/n%
\mathbb{Z}
}(x+y)-\tsum \limits_{x=y\in 
\mathbb{Z}
/n%
\mathbb{Z}
}(x+y)\}=0-0=0.
\end{equation*}

Now, taking the absolute value on both sides of (\ref{CalcDet}), using $%
F_{n}^{4}=n^{2}\cdot Id$ and the positivity of $\tprod
\limits_{j=1}^{n-1}(\sin (\tfrac{\pi \cdot j}{n}))^{n-j},$ gives us 
\begin{equation*}
n^{\frac{n}{2}}=2^{\frac{n(n-1)}{2}}\tprod \limits_{j=1}^{n-1}(\sin (\tfrac{%
\pi \cdot j}{n}))^{n-j},
\end{equation*}%
hence $\det \left( F_{n}\right) =i^{\frac{n(n-1)}{2}}n^{\frac{n}{2}}$. This
concludes the proof of the proposition.

\subsection{Proof of Theorem \protect \ref{DFT_lemma}}

First, we explain why the operator $\rho _{n}\left( \mathrm{w}\right) $ does
not depend on the choice of linearization $\rho _{n}$. Any two linearization
differ by a character $\chi $ of $SL_{2}\left( 
\mathbb{Z}
/n%
\mathbb{Z}
\right) $ therefore it is enough to show that $\chi \left( \mathrm{w}\right)
=1.$ By Proposition \ref{order_prop}, $\chi ^{n}=1$, hence $\chi \left( 
\mathrm{w}\right) =1$, since $\mathrm{w}^{4}=1$ and $\gcd \left( 4,n\right)
=1$.

Next, we explain the relation $F_{n}=C_{n}\cdot \rho _{n}\left( \mathrm{w}%
\right) $. The operator $\rho _{n}\left( \mathrm{w}\right) $ is
characterized up to a unitary scalar by the identity (see formula (\ref%
{Egorov})) $\rho _{n}\left( \mathrm{w}\right) \pi \left( h\right) \rho
_{n}\left( \mathrm{w}\right) ^{-1}=\pi \left( \mathrm{w}\left( h\right)
\right) $ for every $h\in H$. Explicit computation reveals that for every $%
h\in H$, $F_{n}\circ \pi \left( h\right) \circ F_{n}^{-1}=\pi \left( \mathrm{%
w}\left( h\right) \right) $, which implies that 
\begin{equation*}
F_{n}=C_{n}\cdot \rho _{n}\left( \mathrm{w}\right) .
\end{equation*}

Finally we evaluate the proportionality coefficient $C_{n}$. Computing
determinants one obtains $\det \left( F_{n}\right) =C_{n}^{n}\cdot \det
\left( \rho _{n}\left( \mathrm{w}\right) \right) $. Now, by Proposition \ref%
{order_prop}, the character $\chi =\det \circ \rho _{n}$ satisfies $\chi
^{n}=1,$since $\mathrm{w}^{4}=1$ and $\gcd \left( 4,n\right) =1$ it implies
that $\chi \left( \mathrm{w}\right) =1$. Hence 
\begin{equation}
\det \left( F_{n}\right) =C_{n}^{n}.  \label{detC}
\end{equation}

The relations $F_{n}^{4}=n^{2}\cdot Id$ and $\rho _{n}\left( \mathrm{w}%
\right) ^{4}=Id$ imply that $C_{n}^{4}=n^{2}$. By Proposition \ref{detF},
Equation (\ref{detC}) and using $\gcd \left( 4,n\right) =1$ one obtains that 
$C_{n}=i^{\frac{n-1}{2}}\sqrt{n}$.

This concludes the proof of the Theorem.

\subsection{Proof of Lemma \protect \ref{ident_lemma}}

Consider the isomorphism $%
\mathbb{Z}
/p%
\mathbb{Z}
\times 
\mathbb{Z}
/q%
\mathbb{Z}
\overset{\simeq }{\rightarrow }%
\mathbb{Z}
/pq%
\mathbb{Z}
$, given by $\left( x,y\right) \longmapsto x\cdot q+y\cdot p$. Now write 
\begin{eqnarray*}
G_{pq} &=&\tsum \limits_{z\in \mathbb{%
\mathbb{Z}
}/pq%
\mathbb{Z}
}e^{\frac{2\pi i}{pq}z^{2}}=\tsum \limits_{x\in \mathbb{%
\mathbb{Z}
}/p%
\mathbb{Z}
}\tsum \limits_{y\in 
\mathbb{Z}
/q%
\mathbb{Z}
}e^{\frac{2\pi i}{pq}\left( x\cdot q+y\cdot p\right) ^{2}} \\
&=&\tsum \limits_{x\in \mathbb{F}_{p}}e^{\frac{2\pi i}{p}qx^{2}}\tsum%
\limits_{y\in \mathbb{F}_{q}}e^{\frac{2\pi i}{q}py^{2}}=\QOVERD( )
{p}{q}\QOVERD( ) {q}{p}G_{p}\cdot G_{q}\text{, }
\end{eqnarray*}%
where in the last equality we used

\begin{lemma}
\label{tech_lemma}For every $a\in \mathbb{F}_{p}^{\ast }$ 
\begin{equation*}
\tsum \limits_{x\in \mathbb{F}_{p}}e^{\frac{2\pi i}{p}ax^{2}}=\QOVERD( )
{a}{p}\tsum \limits_{x\in \mathbb{F}_{p}}e^{\frac{2\pi i}{p}x^{2}}.
\end{equation*}
\end{lemma}

This concludes the proof of the Lemma \ref{ident_lemma}.

\subsubsection{Proof of Lemma \protect \ref{tech_lemma}}

The statement follows from the following basic identity%
\begin{equation}
\tsum \limits_{x\in \mathbb{F}_{p}}e^{\frac{2\pi i}{p}ax^{2}}=\tsum
\limits_{x\in \mathbb{F}_{p}}e^{\frac{2\pi i}{p}ax}\QOVERD( ) {x}{p}\text{,}
\label{identity}
\end{equation}

which can be explained by observing that both sides are equal to $\tsum
\limits_{x\in \mathbb{F}_{p}}e^{\frac{2\pi i}{p}ax}(1+\QOVERD( ) {x}{p}).$

Now using (\ref{identity}) we can write 
\begin{eqnarray*}
\tsum \limits_{x\in \mathbb{F}_{p}}e^{\frac{2\pi i}{p}ax^{2}} &=&\tsum
\limits_{x\in \mathbb{F}_{p}}e^{\frac{2\pi i}{p}ax}\QOVERD( ) {x}{p}=\tsum
\limits_{x\in \mathbb{F}_{p}}e^{\frac{2\pi i}{p}x}\QOVERD( ) {a^{-1}\cdot
x}{p} \\
&=&\QOVERD( ) {a^{-1}}{p}\tsum \limits_{x\in \mathbb{F}_{p}}e^{\frac{2\pi i}{%
p}x}\QOVERD( ) {x}{p}=\QOVERD( ) {a}{p}\tsum \limits_{x\in \mathbb{F}_{p}}e^{%
\frac{2\pi i}{p}x^{2}},
\end{eqnarray*}%
where, in the second equality we applied a change of variables $x\mapsto ax$%
. This concludes the proof of the lemma.

\subsection{Proof of Proposition \protect \ref{trace_prop}}

By Proposition \ref{tensor_prop}, the representations $\rho _{p}\otimes \rho
_{q}$ and $\rho _{pq}\left[ \psi _{p+q}\right] $ differ by a character $\chi 
$ of the group $SL_{2}\left( 
\mathbb{Z}
/pq%
\mathbb{Z}
\right) $. Since $\chi $ is of odd order (see Proposition \ref{order_prop} )
and $\mathrm{w}^{4}=1$, we have $\chi _{1}\left( \mathrm{w}\right) =1$.
Consequently, we get 
\begin{equation*}
Tr\left( \rho _{pq}\left[ \psi _{p+q}\right] (\mathrm{w})\right) =Tr\left(
\rho _{p}(\mathrm{w})\right) \cdot Tr\left( \rho _{q}(\mathrm{w})\right) .
\end{equation*}

Let $S\in GL_{2}\left( 
\mathbb{Z}
/pq%
\mathbb{Z}
\right) $ such that $\det S=p+q$. Define the conjugate representation $%
Ad_{S}\rho _{pq}$ by%
\begin{equation*}
Ad_{S}\rho _{pq}\left( g\right) =\rho _{pq}\left( Ad_{S}g\right)
\end{equation*}

It is not difficult to show that the representations $\rho _{pq}\left[ \psi
_{p+q}\right] $ and $Ad_{S}\rho _{pq}$ differ by a character $\mu $ of the
group $SL_{2}\left( 
\mathbb{Z}
/pq%
\mathbb{Z}
\right) $. Again, since $\mu $ is of odd order (see Proposition \ref%
{order_prop} ) and $\mathrm{w}^{4}=1$, we have $\mu \left( \mathrm{w}\right)
=1$, which implies that 
\begin{equation*}
Tr\left( \rho _{pq}\left[ \psi _{p+q}\right] (\mathrm{w})\right) =Tr\left(
\rho _{pq}\left( Ad_{S}\mathrm{w}\right) \right) \text{.}
\end{equation*}

It is enough to show that there exists an element $g\in SL_{2}\left( 
\mathbb{Z}
/pq%
\mathbb{Z}
\right) $ such that 
\begin{equation}
Ad_{g}\mathrm{w}=Ad_{S}\mathrm{w}\text{.}  \label{conj_eq}
\end{equation}

This follows from the following general statement about conjugacy classes of
regular semisimple elements in $SL_{2}(\mathbb{%
\mathbb{Z}
}/n%
\mathbb{Z}
)$: Let $n\in 
\mathbb{N}
$ be an odd integer. Consider the natural homomorphism 
\begin{equation*}
SL_{2}(\mathbb{%
\mathbb{Z}
}/n%
\mathbb{Z}
)\rightarrow \prod_{p|n}SL_{2}(\mathbb{F}_{p})\text{.}
\end{equation*}

\begin{definition}
An element $g\in SL_{2}(\mathbb{%
\mathbb{Z}
}/n%
\mathbb{Z}
)$ is called regular semisimple if its image $g_{p}\in SL_{2}(\mathbb{F}%
_{p}) $ is regular semisimple, for every $p|n$.
\end{definition}

\begin{lemma}
\label{semisimple_lemma} Let $g_{0}\in SL_{2}(\mathbb{%
\mathbb{Z}
}/n%
\mathbb{Z}
)$ be a regular semisimple element and $S\in GL_{2}\left( \mathbb{%
\mathbb{Z}
}/n%
\mathbb{Z}
\right) $ then there exists an element $g\in SL_{2}(\mathbb{%
\mathbb{Z}
}/n%
\mathbb{Z}
)$ such that 
\begin{equation*}
Ad_{g}\mathrm{w}=Ad_{S}\mathrm{w}\text{.}
\end{equation*}
\end{lemma}

Invoking Lemma \ref{semisimple_lemma}, Equation (\ref{conj_eq}) now follows
from the fact that $\mathrm{w}$ is a regular semisimple element in $%
SL_{2}\left( 
\mathbb{Z}
/pq%
\mathbb{Z}
\right) $.

This concludes the proof of the proposition.

\subsubsection{Proof of Lemma \protect \ref{semisimple_lemma}}

By the Chinese reminder theorem it is sufficient to prove the assertion in
the case $n=p^{k}$ with $p$ an odd prime number. Moreover, using standard
lifting argument (Hensel's lemma) the statement can be reduced further to
the case where $n=p$.

Denote $g_{1}=Ad_{S}g_{0}$. Consider the set 
\begin{equation*}
O=\left \{ g\in SL_{2}\left( \mathbb{F}_{p}\right) :gg_{0}=g_{1}g\right \} .
\end{equation*}

Our goal is to show that $O\neq \varnothing $.

Let $\overline{\mathbb{F}}_{p}$ denote an algebraic closure of $\mathbb{F}%
_{p}$ and consider the algebraic variety 
\begin{equation*}
\mathbf{O}=\left \{ g\in SL_{2}\left( \overline{\mathbb{F}}_{p}\right)
:gg_{0}=g_{1}g\right \} \text{.}
\end{equation*}

Since $g_{0},g_{1}\in SL_{2}\left( \mathbb{F}_{p}\right) $, the variety $%
\mathbf{O}$ is defined over the finite field $\mathbb{F}_{p}$ and $O$ can be
naturally identified with the set of rational points $O=\mathbf{O}\left( 
\mathbb{F}_{p}\right) $. Since $g_{0}$ is regular semisimple it can be
easily verified that the variety $\mathbf{O}$ is not empty. Moreover, it is
a principal homogenous space over the centralizer subgroup $Z\left(
g_{0}\right) =\left \{ g\in SL_{2}\left( \overline{\mathbb{F}}_{p}\right)
:gg_{0}g^{-1}=g_{0}\right \} $ which is isomorphic to $\mathbb{G}_{m}$. In
particular $\mathbf{O}$ is connected. This implies that $\mathbf{O}\left( 
\mathbb{F}_{p}\right) \neq \varnothing $.

This concludes the proof of the lemma.

\end{document}